\documentclass[12pt]{report}
\usepackage{amssymb,amsmath,makeidx,verbatim}
\newcommand{\R}{\mathbb{R}}

\newcommand{\N}{\mathbb{N}}
\newcommand{\C}{\mathbb{C}}
\newcommand{\Z}{\mathbb{Z}}

\newcommand{\be}{\begin{enumerate}}
\newcommand{\ee}{\end{enumerate}}
\newcommand{\bq}{\begin{eqnarray*}}
\newcommand{\eq}{\end{eqnarray*}}

\begin{document}
%\pagenumbering{roman}
\newcommand{\disp}{\displaystyle}
\thispagestyle{empty}
\begin{center}
\textsc{On Quasi-Symmetric Polynomially Bounded Fr$\acute{\mbox{e}}$chet Algebras\\}
\ \\
\textsc{Olufemi O. Oyadare}\\
\ \\
Department of Mathematics,\\
Obafemi Awolowo University,\\
Ile-Ife, $220005,$ NIGERIA.\\
\text{E-mail: \textit{femi\_oya@yahoo.com}}\\
\end{center}
\begin{quote}
{\bf Abstract.} {\it This paper concerns the notion of a symmetric algebra and its generalization to a quasi-symmetric algebra. We study the structure of these algebras in respect to their \textit{hull-kernel regularity} and existence of some ideals, especially the \textit{hull-minimal ideals}.}
\end{quote}
\ \\

{\bf \S1. Introduction.} It is well-known, from J. Ludwig ($1998$) that every semisimple symmetric polynomially bounded Fr$\acute{\mbox{e}}$chet algebra is hull-kernel regular and has a hull-minimal ideal generated by some elements of the algebra. This result gives a way of verifying the existence of hull-minimal ideals and of computing their basis elements in non-normable algebras of harmonic analysis, as has been shown for the Schwartz algebras of nilpotent and connected semisimple groups in J. Ludwig ($1998$) and O. O. Oyadare ($2016$) respectively.

In this paper we give a generalization of the notion of a symmetric algebra to that of a quasi-symmetric and establish the importance of this generalization by showing that every semisimple quasi-symmetric polynomially bounded Fr$\acute{\mbox{e}}$chet algebra is hull-kernel regular.

The results contained herein form a part of results of the author's thesis (O. O. Oyadare ($2016$)) at the University of Ibadan.
\ \\
\ \\
\ \\
\ \\
$\overline{2010\; \textmd{Mathematics}}$ Subject Classification: $43A85, \;\; 22E30, \;\; 22E46$
\ \\
\ \\

\ \\
{\bf \S2. Preliminaries.} For the connected semisimple Lie group
$G$ with finite center, we denote its Lie algebra by $\mathfrak{g}$
whose \textit{Cartan decomposition} is given as $\mathfrak{g} = \mathfrak{t}\oplus\mathfrak{p}$.
We also denote by $K$ the analytic subgroup of $G$ with Lie
algebra $\mathfrak{t}.$  $K$ is then a maximal compact subgroup of $G$.
Choose a maximal abelian subspace  $\mathfrak{a}$ of $\mathfrak{p}$ with algebraic
dual $\mathfrak{a}^*$ and set $A =\exp \mathfrak{a}.$  For every $\lambda \in \mathfrak{a}^*$ put
$$\disp\mathfrak{g}_{\lambda} = \{X \in \mathfrak{g}: [H, X] =
\lambda(H)X, \forall  H \in \mathfrak{a}\},$$ and call $\lambda$ a restricted
root of $(\mathfrak{g},\mathfrak{a})$ whenever $\mathfrak{g}_{\lambda}\neq\{0\}$.
Denote by $\mathfrak{a}'$ the open subset of $\mathfrak{a}$
where all restricted roots are $\neq 0$,  and call its connected
components the \textit{Weyl chambers}.  Let $\mathfrak{a}^+$ be one of the Weyl
chambers, define the restricted root $\lambda$ positive whenever it
is positive on $\mathfrak{a}^+$ and denote by $\triangle^+$ the set of all
restricted positive roots.  We then have the \textit{Iwasawa
decomposition} $G = KAN$, where $N$ is the analytic subgroup of $G$
corresponding to $\disp \mathfrak{n} = \sum_{\lambda \in \triangle^+} \mathfrak{g}_{\lambda}$,
and the \textit{polar decomposition} $G = K\cdot
cl(A^+)\cdot K,$ with $A^+ = \exp \mathfrak{a}^+,$ and $cl(A^{+})$ denoting the closure of $A^{+}.$  If we set $\disp M = \{k
\in K: Ad(k)H = H$, $H\in \mathfrak{a}\}$ and $\disp M' = \{k
\in K : Ad(k)\mathfrak{a} \subset \mathfrak{a}\}$ and call them the
\textit{centralizer} and \textit{normalizer} of $\mathfrak{a}$ in $K,$ respectively, then;
(i) $M$ and $M'$ are compact and have the same Lie algebra and
(ii) the factor  $\mathfrak{w} = M'/M$ is a finite group called the \textit{Weyl
group}.  $\mathfrak{w}$ acts on $\mathfrak{a}^*_{\C}$ as a group of linear
transformations by the requirement $(s\lambda)(H) =
\lambda(s^{-1}H)$, $H \in \mathfrak{a}$, $s \in \mathfrak{w}$, $\lambda \in
\mathfrak{a}^*_\mathbb{\C}$, the complexification of $\mathfrak{a}^*$.  We then have the
\textit{Bruhat decomposition} $$\disp G = \bigsqcup_{s\in \mathfrak{w}} B m_sB$$ where
$B = MAN$ is a closed subgroup of $G$ and $m_s \in M'$ is the
representative of $s$ (i.e., $s = m_sM$).\\

Some of the most important functions on $G$ are the \textit{spherical
functions} which we now discuss as follows.  A non-zero continuous
function $\varphi$ on $G$ shall be called a \textit{(zonal) spherical
function} whenever $\varphi(e)=1,$ $\varphi \in C(G//K):=\{g\in
C(G)$: $g(k_1 x k_2) = g(x)$, $k_1,k_2 \in K$, $x \in G\}$ and $f*\varphi
= (f*\varphi)(e)\cdot \varphi$ for every $f \in C_c(G//K),$ where $(f \ast g)(x):=\int_{G}f(y)g(y^{-1}x)dy.$  This
leads to the existence of a homomorphism $\lambda :
C_c(G//K)\rightarrow \C$ given as $\lambda(f) = (f*\varphi)(e)$.
This definition is equivalent to the satisfaction of the functional relation $$\disp
\int_K\varphi(xky)dk = \varphi(x)\varphi(y),\;\;\;x,y\in G.$$  It has
been shown by Harish-Chandra [$8.$] that spherical functions on $G$
can be parametrized by members of $\mathfrak{a}^*_{\C}$.  Indeed every
spherical function on $G$ is of the form $\disp
\varphi_{\lambda}(x) = \int_Ke^{(i\lambda-p)H(xk)}dk,\; \lambda
\in \mathfrak{a}^*_{\C},$  $\disp \rho =
\frac{1}{2}\sum_{\lambda\in\triangle^+} m_{\lambda}\cdot\lambda,$ where
$m_{\lambda}=dim (\mathfrak{g}_\lambda),$ and that $\disp \varphi_{\lambda} =
\varphi_{\mu}$ iff $\lambda = s\mu$ for some $s \in \mathfrak{w}$.  Some of
the well-known properties are $\varphi_{-\lambda}(x^{-1}) =
\varphi_{\lambda}(x)$, $\varphi_{-\lambda}(x) =
\bar{\varphi}_{\bar{\lambda}}(x),\; \lambda \in \mathfrak{a}^*_{\C},\;\;
x \in G$, and if $\Omega$ is the \textit{Casimir operator} on $G$ then
$\Omega\varphi_{\lambda} = -(\langle\lambda,\lambda\rangle +
\langle \rho, \rho\rangle)\varphi_{\lambda},$ where $\lambda \in
\mathfrak{a}^*_{\C}$ and $\langle\lambda,\mu\rangle
:=tr(adH_{\lambda} \ adH_{\mu})$ for elements $H_{\lambda}$, $H_{\mu}
\in {\mathfrak{a}}.$ The elements $H_{\lambda}$, $H_{\mu}
\in {\mathfrak{a}}$  are uniquely defined by the requirement that $\lambda
(H)=tr(adH \ adH_{\lambda})$ and $\mu
(H)=tr(adH \ adH_{\mu})$ for every $H \in {\mathfrak{a}}.$ Clearly $\Omega\varphi_0 = 0.$\\

Let $$\disp \varphi_0(x) = \int_{K}\exp(-\rho(H(xk)))dk$$ be denoted
as $\Xi(x)$ and define $\disp\sigma: G \rightarrow \C$ as
$\disp\sigma(x) = \|X\|$ for every $x = k\exp X \in G,\;\; k \in K,\; X
\in \mathfrak{a}$ where $\|\cdot\|$ is a norm on the finite-dimensional
space $\mathfrak{a}.$ These two functions are spherical functions on
$G$ and there exist numbers $c,d$ such that $1 \leq \Xi(a)
e^{\rho(\log a)} \leq c(1+\sigma(a))^d.$ Also there exists $r_0
> 0$ such that $\disp c_0 =: \int_G\Xi(x)^2(1+\sigma(x))^{r_0}dx
< \infty.$ For each
$0 \leq p \leq 2$ define ${\cal C}^p(G)$ to be the set consisting of
functions $f$ in $C^{\infty}(G)$ for which $$\disp \|f\|_{g_1,
g_2;m} :=\sup_G|f(g_1; x ; g_2)|\Xi (x)^{-2/p}(1+\sigma(x))^m <
\infty$$ where $g_1,g_2 \in \mathfrak{U}(\mathfrak{g}_{\C})$, the \textit{universal
enveloping algebra} of $\mathfrak{g}_{\C}$, $m \in \Z^+, x \in G$,
$\disp f(x;g_2) := \left.\frac{d}{dt}\right|_{t=0}f(x\cdot(\exp tg_2))$
and $\disp f(g_1;x) :=\left.\frac{d}{dt}\right|_{t=0}f((\exp
tg_1)\cdot x)$.  We call ${\cal C}^p(G)$ the Schwartz space on $G$
for each $0 < p \leq 2$ and note that ${\cal C}^2(G)$ is the
well-known Harish-Chandra space of rapidly decreasing functions on
$G.$ The inclusions $C^{\infty}_{c}(G) \subset {\cal C}^p(G)
\subset L^p(G)$ hold and with dense images. It also follows that
${\cal C}^p(G) \subseteq {\cal C}^q(G)$ whenever $0 \leq p \leq q
\leq 2.$ Each ${\cal C}^p(G)$ is closed under \textit{involution} and the
\textit{convolution}, $*.$ Indeed ${\cal C}^p(G)$ is a Fr$\acute{e}$chet algebra. We endow ${\cal C}^p(G//K)$
with the relative topology as a subset of ${\cal C}^p(G).$\\

For any measurable function $f$ on $G$ we define the \textit{spherical
transform} $\hat{f}$ as $\disp\hat{f}(\lambda) = \int_G f(x)
\varphi_{-\lambda}(x)dx,$ $\disp \lambda \in \mathfrak{a}^*_{\C}.$ It
is known that for $f,g \in L^1(G)$ we have
\begin{enumerate}
\item [(i)] $(f*g)^{\wedge} = \hat{f}\cdot\hat{g}$ on $ {\mathfrak{F}}^{1}$
whenever $f$ (or $g$) is right - (or left-) $K$-invariant; \item
[(ii)] $\disp (f^*)^{\wedge}(\varphi) =
\overline{\hat{f}(\varphi^*)}, \varphi \in {\mathfrak{F}}^{1}$; hence
$(f^*)^{\wedge} = \overline{\hat{f}}$ on ${\cal P}$: and, if we
define $\disp f^{\#}(g) := \int_{K\times K}f(k_1xk_2)dk_1dk_2,  x\in
G,$ then \item [(iii)] $(f^{\#})^{\wedge}=\hat{f}$ on ${\mathfrak{F}}^{1}$.
\end{enumerate}

\ \\
{\bf \S3. Main Results.} We start with fixing some notions.\\

{\bf 3.1 Definition.} A Fr$\acute{\mbox{e}}$chet algebra $A$ is said to be involutive if there is a map $a \longmapsto a^*$ on $A$ such that $a^{**} = a$, $(a+b)^* = a^*+b^*,\;(\lambda a)^* = \bar{\lambda} a^*,\; \lambda\in \C,$ and $(a\cdot b)^* = b^*\cdot a^*,$ for all $a,b\in A.$

It is clear that, for each $0 < p \leq 2,$ the Fr$\acute{\mbox{e}}$chet algebra $A = {\cal C}^p(G)$ or ${\cal C}^p(G//K)$ is involutive with involution $f \longmapsto f^*$ given as $$f^*(x) = \overline{f(x^{-1})},\; x \in G.$$

Now since the Fr$\acute{\mbox{e}}$chet algebra we are ultimately going to consider in this work is a Schwartz algebra we need to impose a growth condition on the members of the present abstract Fr$\acute{\mbox{e}}$chet algebra which fits into the general behaviour of Schwartz functions. Motivated by the estimates of members of ${\cal C}^p(G)$ as known above we give the following definition.\\

{\bf 3.2 Definition.} Let $A = (A, \{p_k\})$ be an involutive algebra, where $\{p_k\}$ is the defining collection of seminorms on $A,$ and set $\disp e(a) = \sum^{\infty}_{k=1}\frac{a^k}{k!}$ for $a \in A.$ An element $b \in A$ is said to be \textit{polynomially bounded} if for every $k \in \N$ there is a constant $c_k = c_k(b)>0$ such that
$$p_k(e(i\lambda b)) \leq c_k(1+|\lambda|)^{c_k},\;\; \mbox{for all}\;\; \lambda \in \R.$$

{\bf 3.3 Remarks.}

$(1.)$ This definition may be compared with the weak inequality of the last chapter using the fact $\Xi(x) \leq 1,$ for all $x \in G.$

$(2.)$ Clearly $e(a) = \exp(a)-1.$ See J. Dixmier ($1960$).

As a guiding example we consider the special (but important) case of when $B$ is a Banach algebra, $(B, \|\cdot\|).$ In this case, if $B(b)$ is a maximal abelian closed subalgebra of $B$ containing $b,$ and $\chi$ is a character on $B(b)$ for which $\chi(b) = \mu \in spec(b) (:=$ spectrum of $b$), we know that in general $\mu \in \C.\;i.e.,\;\mu = \alpha+i\beta$ where $\alpha,\beta \in \R.$ However we also know that, for every $\lambda \in \R,$ $\disp |e^{i\lambda\mu}| = |e^{i\lambda(\alpha+i\beta)}| = |e^{i\lambda\alpha}||e^{-\lambda\beta}| = 1|e^{-\lambda\beta}| = e^{-\lambda\beta}$; so that $\disp e^{-\lambda\beta} = |e^{i\lambda\mu}| = |\exp(i\lambda\chi(b)| = |1+e(i\lambda\chi(b))| = |1+\chi(e(i\lambda b))| \leq 1+\|e(i\lambda b)\|.\;i.e.,\;\disp e^{-\lambda\beta} \leq 1+\|e(i\lambda b)\|$ for all $b \in B,$ meaning that $\|e(i\lambda b)\|$ grows exponentially in $\lambda,$ with $\beta$ a real constant.

Thus in order to have an element $b \in B = (B, \|\cdot\|)$ to be polynomially bounded it must be such that $\beta$ in the above inequality must be zero. $i.e.,$ the element $b$ must have real spectrum. Thus, since every Banach algebra is a Fr$\acute{\mbox{e}}$chet algebra we consider the requirement of having a real spectrum for a polynomially bounded element of the Fr$\acute{\mbox{e}}$chet algebra, $(A,\{p_k\}).$ We however recall that an involutive Banach algebra in which the spectrum of every self-adjoint element is a subset of $\R$ is called \textit{symmetric} and then introduce the following notion of symmetricity.\\

{\bf 3.4 Definition.} A Fr$\acute{\mbox{e}}$chet algebra $A$ is said to be \textit{symmetric} if it admits a continuous involution and there exists a continuous *-homomorphism, $\sigma,$ of $A$ into a $C^*-$algebra, $\mathfrak{C},$ such that $spec_A(a) = spec_{\mathfrak{C}}(\sigma(a)),$ for every $a \in A$ (Here $\text{spec}_{A}(a)$ represents the spectrum of an element $a$ in $A$ defined as $\{\lambda \in \C : a- \lambda \cdot 1\;\mbox{is not invertible} \}$).\\

{\bf 3.5 Remarks.}

$(1.)$  By a continuous *-homomorphism $\sigma:A\rightarrow \mathfrak{C}$ we mean a continuous homomorphism $\sigma$ in which $\sigma(a^*) = \sigma(a)^{\#}$ where $A$ and $\mathfrak{C}$ are $*-$ and $\#-$involutive algebras respectively.

$(2.)$ The requirement on the spectrum in $(3.4)$ above is equivalent to saying that the continuous *-homomorphism, $\sigma,$ is spectrum invariant.

$(3.)$ If in $(3.4)$ above, we require only that $spec_A(a) \subseteq spec_{\mathfrak{C}}(\sigma(a)),$ then we shall refer to $A$ a \textit{quasi-symmetric} Fr$\acute{\mbox{e}}$chet algebra. Clearly every symmetric Fr$\acute{\mbox{e}}$chet algebra is automatically quasi-symmetric but not conversely. Thus the notion of quasi-symmetricity of a Fr$\acute{\mbox{e}}$chet algebra is more general than that of the symmetricity in $(3.4).$ Indeed we shall in a moment extend the results of J. Ludwig ($1998:$ Propositions $1.8$ and $1.10$) to include all quasi-symmetric Fr$\acute{\mbox{e}}$chet algebras.

We now consider the notion of a \textit{functional calculus} (E. Hille and R. S. Phillips ($1957$)) that is needed shortly. We start with a Banach algebra example.\\

{\bf 3.6  Definition.} Let $B$ be a Banach *-algebra. A function $\varphi$ is said to operate on an element $f \in B$ if
\begin{enumerate}
\item [$(i.)$] the Gelfand transform, $\hat{f},$ of $f,$ with respect to the smallest commutative Banach sub-algebra $B(f)$ containing $f,$ is real, and
\item [$(ii.)$] there exists a $g\in B(f)$ such that $\varphi \circ \hat{f} = \hat{g}.$
\end{enumerate}

We shall in this case write $\varphi\{f\} = g.$ Indeed, if for $f\in B$ we write $\disp e(f) = \sum^{\infty}_{k=1}\frac{(if)^k}{k!}$ and if $\disp \|e(nf)\|_B = O(\|n\|^N)$ as $\|n\|\rightarrow \infty,$ then every $\varphi \in C^k_c(\R)$ with $k>N+1$, and $\varphi(0)=0$ operates on $f$ and $\disp \varphi\{f\} = \frac{1}{2\pi i}\int_{\R}\hat{\varphi}(\lambda)e(\lambda f)d\lambda$ so that $\disp \|\varphi\{f\}\|_B \leq c\|\varphi\|_{C^k_c(\R)}.$

In a Fr$\acute{\mbox{e}}$chet algebra $A$ we may also employ the functional calculus of $C^{\infty}$-functions on polynomially bounded elements of $A$ ($cf.$ J. Dixmier $(1960)$). Let us denote by $C^{\infty}_{c,0}(\R)$ members $\varphi \in C^{\infty}_c(\R)$ in which $\varphi(0) = 0,$ then the integral $\disp \frac{1}{2\pi i}\int_{\R}\hat{\varphi}(\lambda)e(i\lambda a)d\lambda$ exists and converges in $A,$ for any polynomially bounded elements $a \in A$ (J. Dixmier, $1960:$ p. $18$). We then define $$\disp\varphi\{a\} = \frac{1}{2\pi i}\int_{\R}\hat{\varphi}(\lambda)e(i\lambda a)d\lambda.$$ This functional calculus on polynomially bounded elements of $A$ has the following interesting properties that makes it what we actually need.

Let $a$ be a polynomially bounded element of the Fr$\acute{\mbox{e}}$chet algebra $A,$ and let $A(a)$ be a maximal abelian closed subalgebra containing $a.$ If $\chi$ is any character on $A(a)$ then $\chi(\varphi\{a\}) = \varphi\{\chi(a)\}.$ We thus have
$$\chi((\psi\cdot\varphi)\{a\}) = (\psi\cdot\varphi)\{\chi(a)\} = \psi\{\chi(a)\}\cdot\varphi\{\chi(a)\} = \chi(\psi\{a\})\cdot\chi(\varphi\{a\})$$
for $\psi, \varphi \in C^{\infty}_{c,0}(\R).$ If $A(a)$ is now semisimple. That is, if $\ker(\chi) = \{0\}$ for every character $\chi$ of $A(a),$ then $$\chi((\psi\cdot\varphi)\{a\}) = \chi(\psi\{a\})\cdot\chi(\varphi\{a\}) = \chi(\psi\{a\}\cdot\varphi\{a\})$$ implying $\disp \chi((\psi\cdot\varphi)\{a\} - \psi\{a\}\cdot\varphi\{a\}) = 0$ so that
$$\psi\{a\} \cdot\varphi\{a\} = (\psi\cdot\varphi)\{a\} \cdots \cdots(I).$$

The use to which this calculus is put is contained in the following.\\

{\bf 3.7 Proposition} (J. Ludwig, $1998:$ p. $80$). Let $a$ be a polynomially bounded element of $A.$ Then there exist $\psi,\varphi \in C^{\infty}_{c,0}(\R)$ such that $$\psi\{a\}\cdot\varphi\{a\} = \varphi\{a\}.$$

{\bf Proof.} Since $C^{\infty}_{c,0}(\R)$ is (completely) regular we know that there exist $\psi, \varphi \in C^{\infty}_{c,0}(\R)$ such that $\psi\cdot\varphi=\varphi.$ Using this relation on the right-hand side of (I) above gives the result$.\;\;\;\Box$

We shall soon see how the result of Proposition $3.7$ above simplifies matters in the proof of an important result of this section, thus making it central to our discussion. In an attempt to generalise this calculus to the algebra ${\cal C}^p(G//K)$ one may introduce a \textit{distributional calculus} on members of ${\cal C}^p(G),\; 0 < p \leq 2.$

This generalisation is analogous to the generalisation of characters of finite and compact groups as functions on the groups to global characters on connected semisimple Lie groups as distributions in which J. G. Arthur ($1974$), W. H. Barker ($1975,$ $1976$ and $1984$) and O. O. Oyadare ($2015$) may prove very useful. However due to the well-developed theory of $\mathfrak{Z}(\mathfrak{g}_{\C})-$finite $K-$finite functions and of cusp forms on connected semisimple Lie groups no calculus or special Fourier transforms is needed in the construction of a basis for $j(C)$ in the Schwartz algebras of focus in the last section of this chapter.

We have seen how to express a primitive Ideal of an associative algebra, $A,$ as the kernel of some algebraically irreducible representation of $A$ ($cf.\;(4.2.1)$). Now that we have the algebra $A$ to be Fr$\acute{\mbox{e}}$chet, on which we have a topology induced by its collection of seminorms, we may employ both algebraic and topological irreducibility of a representation of $A$ and make a comparison between them. Indeed we have the following generalisation of $(1.8)$ in J. Ludwig ($1998:$ p. $18$) which is needed in the proof of the major result of this section.\\

{\bf 3.8 Lemma.} Every algebraically irreducible representation space of a quasi-symmetric Fr$\acute{\mbox{e}}$chet algebra $A$ is equivalent to a submodule of a topologically irreducible representation space of $A.$

{\bf Proof.} Since $A$ is quasi-symmetric so also is the adjunction $\tilde{A} :=\C 1\oplus A,$ of 1 to $A.$ Thus we may assume, without any loss of generality, that $A$ and $\mathfrak{C}$ have identities. Now let $\sigma$ be as in $(3.4)$ and let $(T,V)$ be an algebraically irreducible representation of $A.$ We claim that $\ker(\sigma)\subseteq \ker(T).$

Indeed, if $x \in \ker(\sigma)$ and $y\in A,$ then the spectrum of $yx$ in $A$ is reduced to $\{0\}.$ Suppose on the contrary to the claim, that $x \in \ker(\sigma)$ does not imply $x \in \ker(T),$ then $Tx \neq 0;$ so that there exists $0 \neq v \in V$ such that $(Tx)v \neq 0;$ and since $T$ is simple, being algebraically irreducible, then we can find an element $y \in A$ such that $T(y)(T(x)v)=v.\;i.e.,\;(T(yx)-\lambda\cdot 1)v = 0.$ This means that 1 is in the spectrum of $yx,$ a contradiction to the fact that the spectrum of $yx$ is reduced to $\{0\}.$ Hence $\ker(\sigma)\subseteq\ker(T).$

It then follows that there exists a proper maximal left Ideal $M$ of $A$ in which $\ker(T) \subseteq M$ in which the simple module $(T, V)$ is equivalent to the left-regular representation of $A$ on $A/M.$ We note that the sum of $\C 1$ and $\sigma(M)$ is direct in $\mathfrak{C},$ since otherwise $1 \in M\bmod(\ker\sigma)$ which implies $1\in\ker\sigma\subseteq(T)\subseteq M.\;i.e.,\;1 \in M,$ which is impossible (see C. E. Rickart, $1974:$ Corollary $2.1.2$).

Define $\tilde{M} := \sigma(\C 1+M)\; (=\C 1 + \sigma(M)\subseteq \mathfrak{C})$ and define a linear functional $\varphi$ on $\tilde{M}$ by the requirement that $\varphi(\lambda\cdot 1+\sigma(m))=\lambda$ for every $\lambda \in \C,\;m \in M.$ In other words for every $x \in M$, which is of the form $x = \lambda\cdot 1+m$ we define $\varphi(\sigma(x))=\lambda.$ Now since $x = \lambda\cdot 1+m \in M$ we have $x-\lambda\cdot 1 = m \in M;$ meaning that $(x-\lambda\cdot 1)$ is non-invertible in $A.\;i.e.,\;\lambda \in spec_A(x),$ and by the hypothesis of quasi-symmetricity on $A,$ we conclude that $\lambda \in Spec_{\mathfrak{C}}(\sigma(x)).$ Thus $|\varphi(\sigma(x))| = |\lambda| \leq \sup\{|\mu| : \mu \in spec_{\mathfrak{C}}(\sigma(x))\} = \|\sigma(x)\|_{\mathfrak{C}}$ ($cf.$ Theorem $8$ of F. F. Bonsall and J. Duncan, $1973:$ p. $23$). Therefore
$$\|\varphi\|_{op} =  \sup_{\begin{array}{c} x\in M\\ x\neq 0 \end{array}} \frac{|\varphi(\sigma(x))|}{\|\sigma(x)\|_{\mathfrak{C}}} \leq 1.$$
Hence by the Hahn-Banach theorem there exists a continuous extension, say $\tilde{\varphi},$ of $\varphi$ to the whole of $\mathfrak{C}$ of norm $\leq 1.$ Since $\varphi(1) = 1$ and $\|\tilde{\varphi}\|_{op}\leq 1,$ then $\tilde{\varphi}$ is a positive linear functional for which $\tilde{\varphi}(\sigma(\C 1+M)) = \tilde{\varphi}(\tilde{M}) = \{0\}.$ Since $M$ is maximal, we then have that $\varphi(\sigma(M)) = \{0\}$ and we can deduce that $M = \{y \in A :\tilde{\varphi}(\sigma(y^*y))=0\}$ which, in particular, shows that $M$ is closed. (Indeed, every proper maximal Ideal is closed).

We therefore put a pre-Hilbert structure, $\langle \cdot,\cdot\rangle$ on $A/M$ by setting $\langle x+M, y+M\rangle :=\tilde{\varphi}(\sigma(y^*x))$.  Let ${\cal H}$ be the completion of $A/M = (A/M, \langle \cdot, \cdot\rangle),$ then the above left-regular representation of $A$ on $A/M$ extends to a unitary representation $\pi$ of $A$ on ${\cal H}$ ($cf.$ J. Dixmier, $1964:$ $(2.4.4)$). Since we may assume that $\tilde{\varphi}$ is a \textit{pure state} (see R. V. Kadison and J. R. Ringrose, $1983:$ p. $213$), we also know that $\pi$ is topologically irreducible ($cf.$ R. V. Kadison and J. R. Ringrose, $1983:$ $(2.5)$)$.\;\;\;\Box$

It follows from the above result that a quasi-symmetric Fr$\acute{\mbox{e}}$chet algebra $A$ has sufficiently many algebraically irreducible (unitary) representations. This allows the use $(2.2.9)\;(i.)$ of C. E. Rickart $(1974)$ to define members of $Prim(A)$ which leads to the discussion of the hull-minimal Ideals.

In order to then state the major result of this section we put the notion of $(3.2)$ in the proper form in which it is needed.\\

{\bf 3.9 Definition.} An involutive Fr$\acute{\mbox{e}}$chet algebra $A$ is said to be \textit{polynomially bounded} if the set $A_0$ of self-adjoint polynomially bounded elements of $A$ is dense in the real subspace $A_n$ of hermitian elements of $A.$

We note here that an involutive Schwartz algebra is automatically a polynomially bounded Fr$\acute{\mbox{e}}$chet algebra. In particular the family of algebras given in $(3.4)$ as ${\cal C}^p(G//K)$ is a family of polynomially bounded Fr$\acute{\mbox{e}}$chet  algebras. We now give a generalisation of $(1.10)$ in J. Ludwig ($1998:$ p. $81$), to a semisimple quasi-symmetric polynomially bounded Fr$\acute{\mbox{e}}$chet algebra. This result is the first major result of this section of the work. It assures us that there is a reasonable framework for the study of the hull-minimal Ideals in a Fr$\acute{\mbox{e}}$chet algebra.\\

{\bf 3.10 Theorem.} Every semisimple quasi-symmetric polynomially bounded Fr$\acute{\mbox{e}}$chet algebra is hull-kernel regular.

{\bf Proof.} Let $A$ be a semisimple quasi-symmetric polynomially bounded Fr$\acute{\mbox{e}}$chet algebra, then by $(3.8)$ for every $J \in Prim(A)$ (which is of the form $J = $ kernel of an algebraically irreducible representation of $A$) there exists a topologically irreducible unitary representation $(\pi_J, {\cal H}_J)$ of $A$ such that $\ker(\pi_J) = J.$

Now let $C$ be any closed subset of $Prim(A)$ and fix $J \in Prim(A)\backslash C$ we need to show that there are $a_J, b_J \in A$ in which
\begin{enumerate}
\item [$(i.)$] $\disp b_J \in \bigcap_{J'\in C}J' (=: \ker(C)),\;\; a_J \notin J,$ and
\item [$(ii.)$] $\disp b_{J}\cdot a_J = a_J.$
\end{enumerate}

Since $C$ is closed there exists $u \in A$ such that $u \in \ker(C)\;\left(=\bigcap_{J'\in C}J'\right)$ and $u \not\in J.$ Equivalently, $\pi_{J'}(u^*u) = 0$ for every $J'\in C$ and $\pi_J(u^*u)\neq 0,$ respectively. By unitarity of $\pi_{J'}$ it is possible, after multiplying $v := u^*u$ with a positive constant, to have $\|\pi_J(v)\|_{op}=1.$ Clearly $v \in A_h$ as $v^* = v.$ Also by the continuity of the *-homomorphism $\sigma$ (as in $(3.4)$) we have that there exist a continuous seminorm $p$ on $A$ such that $\|\sigma(a)\|_{\mathfrak{C}} \leq p(a),$ for all $a \in A.$ Thus for any unitary representation $\pi$ of $A$ we have $\|\pi(a)\|_{op} \leq \|\sigma(a)\|_{\mathfrak{C}} \leq p(a),\; a \in A.$

As $A_0$ is dense in $A_n$ we can choose $a_0 \in A_0$ such that $p(a_0-v) < 1/m,\;\; m = 2,3,4,\dots,$ and real $C^{\infty}$-functions $\varphi, \psi$ such that $\psi$ vanishes in a neighbourhood ${\cal N}_\psi$ of $\disp \left[-\frac{2}{m},\; \frac{2}{m}\right], \;\varphi=1$ on $\disp \left[1-\frac{1}{m},\; 1+\frac{1}{m}\right]$ and $\psi\cdot\varphi = \varphi.$ If we now set $\disp b_J = \psi\{a_0\}$ and $a_J = \varphi\{a_0\}$, we have, by the functional calculus on semisimple $A,$ that $\disp \psi\{a_0\}\cdot\varphi\{a_0\} = \varphi\{a_0\}$\ \ \ \ \ ($cf.$ $(3.7)$). $i.e.,$ $b_J \cdot a_J = a_J,$ as required in $(ii.)$ above.

To now verify $(i.)$ we note that
\begin{eqnarray*}
\|\pi_{J'}(a_0)\|_{op} &=& \|\pi_{J'}(a_0-v)\|_{op}\;\; \mbox{(since $\disp \pi_{J'}(v) = \pi_{J'}(u^*u) = 0,$ for every $J'\in C$)}\\
& \leq & p(a_0-v) \\
& \leq & \frac{1}{m},\;\;
\end{eqnarray*}
which shows clearly that $\pi_{J'}(a_0) \in {\cal N}_\psi.$ Thus
$$\psi(\pi_{J'}(a_0)) = 0.\;i.e.,\;\disp \pi_{J'}(\psi\{a_0\}) = 0$$ implying that $\pi_{J'}(b_J)=0.\;i.e.,\;\disp b_J \in \ker(\pi_{J'}) = J'.$ In short, $b_J \in J'$ for each $J' \in C.\;i.e.,\;\disp b_J \in \bigcap_{J'\in C}J' = \ker(C),$ as required. Also
\begin{eqnarray*}
|\|\pi_J(a_0)\|_{op}-1| &=& |\|\pi_J(a_0)\|_{op}-\|\pi_J(v)\|_{op}|\\
& \leq & \|\pi_J(a_0)-\pi_J(v)\|_{op}\;\;\; \mbox{(by continuity of $\|\cdot\|_{op}$)}\\
&=& \|\pi_J(a_0-v)\|_{op} \leq p(a_0-v) \leq 1/m,\;\; m = 2,3,4,\dots;
\end{eqnarray*}
implying that $\disp 1-\frac{1}{m} \leq \|\pi_J(a_0)\|_{op} \leq 1+\frac{1}{m}.\;i.e.,\;\disp\varphi(\pi_J(a_0))\neq 0.$ Thus $\disp \pi_J(a_J) = \pi_J(\varphi\{a_0\}) = \varphi(\pi_J(a_0)) \neq 0.\;i.e.,\;a_J \notin\ker(\pi_J)=J.\;\;\;\Box$\\

{\bf 3.11 Corollary.} Let $A$ be a semisimple quasi-symmetric polynomially bounded Fr$\acute{\mbox{e}}$chet algebra and let $C$ be a closed subset of $Prim(A).$ Then the hull-minimal Ideal, $j(C),$ exists and is generated by the elements $a_J,\; J \notin C.$

{\bf Proof.} The proof follows if we combine Lemma $3.11$ with Theorem $2.3$ contained in U. N. Bassey and O. O. Oyadare ($2011$) (or with Theorem $4.2.20$ of O. O. Oyadare (2016))$.\;\;\;\Box$\\

Lemma $3.8,$ Theorem $3.10$ and Corollary $3.11$ hold in particular for the Harish-Chandra-type algebras $A=\mathcal{C}^{p}(G).$

The last result is an improvement on $(1.10)$ of J. Ludwig $(1998)$ as afforded by the general notion of quasi-symmetricity introduced in $(4.3.5)(3.)$ above. We refer to U. N. Bassey and O. O. Oyadare($2011$), as well as O. O. Oyadare ($2016$), for the explicit structure of the basis elements, $a_{J},$ of (type-I) hull-minimal ideals in the quasi-symmetric polynomially bounded Fr$\acute{\mbox{e}}$chet algebras, $\mathcal{C}^{p}(G),$ of a (connected) semisimple $G.$

\ \\

{\bf Acknowledgement}: The author wishes to appreciate his Ph. D thesis adviser, Dr. U. N. Bassey (University of Ibadan, Nigeria).\\
\ \\

{\bf   References.}
\begin{description}
\item [{[1.]}] Arthur, J. G., Some tempered distributions on semisimple groups of real rank one. \textit{Ann. of Math.} 100  (1974): $553$-$584.$
        \item [{[2.]}] Barker, W. H., The spherical Bochner theorem on semisimple Lie groups. \textit{J. Funct. Anal.} 20 (1975): $179$-$207.$
            \item [{[3.]}] Barker, W. H., Positive definite distributions on unimodular Lie groups. \textit{Duke Math. J.} 43. 1 (1976): $71$-$79,$
                \item [{[4.]}] Barker, W. H., Tempered, invariant, positive-definite distributions on $SU(1,1)/\{\pm 1\}.$ \textit{Illinois J. Math.} 28. 1 (1984): $83$-$102.$
                    \item [{[5.]}] Bassey, U. N. and Oyadare, O. O., On Hull-minimal ideals in the Schwartz algebras of spherical functions on a semisimple Lie group. \textit{International Journal of Functional Analysis, Operator Theory and Applications.} $3\;(1)\;(2011):$ $37-52.$
                    \item [{[6.]}] Bonsall, F. F. and Duncan, J., \textit{Complete normed algebras.} Berlin: Springer-Verlag, $1973.$
                    \item [{[7.]}] Dixmier, J., Op$\acute{\mbox{e}}$rateurs de rang fini dans les repr$\acute{\mbox{e}}$sentations unitaires. \textit{Inst. Hautes $\acute{\mbox{E}}$tudes.  Sci. Publ. Math.} 6 (1960): $13$-$25.$
\item [{[8.]}] Dixmier, J., \textit{Les C*-algebres et leurs representations.} Paris: Gauthier-Villars, 1964.
                \item [{[9.]}] Harish-Chandra, Spherical functions on a semisimple Lie group, I. \textit{Amer. J. Math.} 80. 2 (1958): $241$-$310.$
                        \item [{[10.]}] Helgason, S., \textit{Differential geometry, Lie groups and symmetric spaces.} New York: Academic Press, 1978.
                            \item [{[11.]}] Hille, E. and Phillips, R. S., \textit{Functional analysis and semi-groups.} {\bf 31.} Rhode Island: Amer. Math. Soc., Colloq. Publ., Providence, 1957.
                                \item [{[12.]}] Kadison, R. V. and Ringrose, J. R., \textit{Fundamentals of the theory of operators algebras,} $vol.$ {\bf 1}, Elementary theory. London: Academic Press, 1983.
                                    \item [{[13.]}] Ludwig, J., Hull-minimal Ideals in the Schwartz algebra of the Heisenberg group, \textit{Studia Mathematica,} {\bf 130}, 1 (1998), pp. $77$-$98.$
                                        \item [{[14.]}] Oyadare, O. O., On harmonic analysis of spherical convolutions on semisimple Lie groups, \textit{Theoretical Mathematics and Applications,} \textbf{5} (\textbf{3}) (2015), p. $19-36.$
                                             \item [{[15.]}] Oyadare, O. O., \textit{On ideal theory for the Schwartz algebras of spherical functions on semisimple Lie groups.} Ph.D Thesis, University of Ibadan (2016). $x+383$pp.
                                                 \item [{[16.]}] Rickart, C. E., \textit{General theory of Banach algebras.} New York: Van Nostrand, 1974.
                        \end{description}
\end{document}